\documentclass{amsart}%
\usepackage{amssymb,amsmath,amsfonts,latexsym,amsthm,geometry,graphicx}
\usepackage{amsmath}
\usepackage{amsfonts}
\usepackage{amssymb}
\usepackage{graphicx}%
\providecommand{\U}[1]{\protect\rule{.1in}{.1in}}
\providecommand{\U}[1]{\protect\rule{.1in}{.1in}}
\providecommand{\U}[1]{\protect\rule{.1in}{.1in}}
\newtheorem{theorem}{Theorem}[section]

\theoremstyle{definition}

\begin{document}
\title[Lower bounds for the constants of the Hardy-Littlewood inequalities]{Lower bounds for the constants of the Hardy-Littlewood inequalities}
\date{}
\author[G. Ara\'{u}jo]{Gustavo Ara\'{u}jo}
\address{Departamento de Matem\'{a}tica \\
Universidade Federal da Para\'{\i}ba \\
58.051-900 - Jo\~{a}o Pessoa, Brazil.}
\email{gdasaraujo@gmail.com}
\author[D. Pellegrino]{Daniel Pellegrino}
\address{Departamento de Matem\'{a}tica \\
Universidade Federal da Para\'{\i}ba \\
58.051-900 - Jo\~{a}o Pessoa, Brazil.}
\email{pellegrino@pq.cnpq.br and dmpellegrino@gmail.com}
\thanks{The authors are supported by CNPq Grant 313797/2013-7 - PVE - Linha 2}

\begin{abstract}
Given an integer $m\geq2$, the Hardy--Littlewood inequality (for real scalars)
says that for all $2m\leq p\leq\infty$, there exists a constant $C_{m,p}%
^{\mathbb{R}}\geq1$ such that, for all continuous $m$--linear forms
$A:\ell_{p}^{N}\times\cdots\times\ell_{p}^{N}\rightarrow\mathbb{R}$ and all
positive integers $N$,
\[
\left(  \sum_{j_{1},...,j_{m}=1}^{N}\left\vert A(e_{j_{1}},...,e_{j_{m}%
})\right\vert ^{\frac{2mp}{mp+p-2m}}\right)  ^{\frac{mp+p-2m}{2mp}}\leq
C_{m,p}^{\mathbb{R}}\left\Vert A\right\Vert .
\]
The limiting case $p=\infty$ is the well-known Bohnenblust--Hille inequality;
the behavior of the constants $C_{m,p}^{\mathbb{R}}$ is an open problem. In
this note we provide nontrivial lower bounds for these constants.

\end{abstract}
\maketitle

\section{Introduction}

Let $\mathbb{K}$ denote the field of real or complex scalars. The multilinear
Bohnenblust--Hille inequality asserts that for all positive integers $m\geq2$
there exists a constant $C_{m}^{\mathbb{K}}\geq1$ such that
\[
\left(  \sum_{j_{1},...,j_{m}=1}^{N}\left\vert A(e_{j_{1}},...,e_{j_{m}%
})\right\vert ^{\frac{2m}{m+1}}\right)  ^{\frac{m+1}{2m}}\leq C_{m}%
^{\mathbb{K}}\Vert A\Vert
\]
for all continuous $m$--linear forms $A:\ell_{\infty}^{N}\times\cdots
\times\ell_{\infty}^{N}\rightarrow\mathbb{K}$ and all positive integers $N$.
The precise growth of the constants $C_{m}^{\mathbb{K}}$ is important for
applications (see \cite{montanaro}) and remains a big open problem. Only very
recently, in \cite{bohr,npss} it was shown that the constants have a
subpolynomial growth. For real scalars, in 2014 (see \cite{diniz}) it was
shown that the optimal constant for $m=2$ is $\sqrt{2}$ and in general
$C_{m}^{\mathbb{R}}\geq2^{1-\frac{1}{m}}$. In the case of complex scalars it
is still an open problem whether the optimal constants are strictly grater
than $1$; in the polynomial case, in 2013 D. N\'{u}\~{n}ez-Alarc\'{o}n proved
that the complex constants are strictly greater than $1$ (see \cite{nunez}).

Even basic questions related to the constants $C_{m}^{\mathbb{K}}$ remain
unsolved. For instance:

\begin{itemize}
\item Is the sequence of optimal constants $\left(  C_{m}^{\mathbb{K}}\right)
_{m=1}^{\infty}$ increasing?

\item Is the sequence of optimal constants $\left(  C_{m}^{\mathbb{K}}\right)
_{m=1}^{\infty}$ bounded?

\item Is $C_{m}^{\mathbb{C}}=1$?
\end{itemize}

The Hardy-Littlewood inequalities are a generalization of the
Bohnenblust--Hille inequality to $\ell_{p}$ spaces. The bilinear case was
proved by Hardy and Littlewood in 1934 (see \cite{hardy}) and in 1981 it was
extended to multilinear operators by Praciano--Pereira (see \cite{pra}).

\begin{theorem}
[Hardy and Littlewood inequalities (\cite{hardy, pra})]\label{000}Let $m\geq2$
be a positive integer and $p\geq2m$. Then there exists a constant
$C_{m,p}^{\mathbb{K}}\geq1$ such that, for every continuous $m$--linear form
$A:\ell_{p}^{N}\times\cdots\times\ell_{p}^{N}\rightarrow\mathbb{K}$,
\[
\left(  \sum_{j_{1},...,j_{m}=1}^{N}\left\vert A(e_{j_{1}},...,e_{j_{m}%
})\right\vert ^{\frac{2mp}{mp+p-2m}}\right)  ^{\frac{mp+p-2m}{2mp}}\leq
C_{m,p}^{\mathbb{K}}\left\Vert A\right\Vert
\]
for all positive integers $N$.
\end{theorem}

It was recently shown that $\left(  C_{m,p}^{\mathbb{K}}\right)
_{m=1}^{\infty}$ is sublinear for $p\geq m^{2}$. More precisely, it was shown
that (see \cite[page 1887]{aps})
\[%
\begin{array}
[c]{llll}%
C_{m,p}^{\mathbb{C}} & \leq\left(  \frac{2}{\sqrt{\pi}}\right)  ^{\frac
{2m(m-1)}{p}}\left(  \displaystyle\prod\limits_{j=2}^{m}\Gamma\left(
2-\frac{1}{j}\right)  ^{\frac{j}{2-2j}}\right)  ^{\frac{p-2m}{p}}, &
\vspace{0.2cm} & \\
C_{m,p}^{\mathbb{R}} & \leq\left(  \sqrt{2}\right)  ^{\frac{2m\left(
m-1\right)  }{p}}\left(  2^{\frac{446381}{55440}-\frac{m}{2}}%
\displaystyle\prod\limits_{j=14}^{m}\left(  \frac{\Gamma\left(  \frac{3}%
{2}-\frac{1}{j}\right)  }{\sqrt{\pi}}\right)  ^{\frac{j}{2-2j}}\right)
^{\frac{p-2m}{p}}, & \text{ for }m\geq14,\vspace{0.2cm} & \\
C_{m,p}^{\mathbb{R}} & \leq\left(  \sqrt{2}\right)  ^{\frac{2m\left(
m-1\right)  }{p}}\left(  \displaystyle\prod\limits_{j=2}^{m}2^{\frac{1}{2j-2}%
}\right)  ^{\frac{p-2m}{p}}, & \text{ for }2\leq m\leq13. &
\end{array}
\]

The precise estimates of the constants of the Hardy--Littlewood inequalities
are unknown and even its asymptotic growth is a mystery (as it happens with
the Bohnenblust--Hille inequality). In this note we provide nontrivial lower
bounds for these inequalities. Following the lines of \cite{diniz}, it is
possible to prove that
\begin{equation}
C_{m,p}^{\mathbb{R}}\geq2^{\frac{mp+2m-2m^{2}-p}{mp}}>1 \label{001}%
\end{equation}
when $p>2m$, but note that when $p=2m$ we have $2^{\frac{mp+2m-2m^{2}-p}{mp}%
}=1$ and thus we do not have nontrivial information. In this paper we also
treat the extreme case $p=2m.$

\begin{theorem}
\label{pop}The optimal constants of the Hardy-Littlewood inequalities satisfy
\[
C_{m,p}^{\mathbb{R}}>1
\]
%whenever $p=2m$ with $p\leq100$ and for all $p>2m.$

\end{theorem}

\section{The proof of Theorem \ref{pop}}

All that it left to prove is the case $p=2m$. We divide the proof in five steps.

\bigskip

\noindent\textbf{Step 1.} Induction.

\bigskip

This first step follows the lines of \cite{diniz}. For $2m\leq p\leq\infty$,
consider
\[%
\begin{array}
[c]{ccccl}%
T_{2,p} & : & \ell_{p}^{2}\times\ell_{p}^{2} & \rightarrow & \mathbb{R}\\
&  & (x^{(1)},x^{(2)}) & \mapsto & x_{1}^{(1)}x_{1}^{(2)}+x_{1}^{(1)}%
x_{2}^{(2)}+x_{2}^{(1)}x_{1}^{(2)}-x_{2}^{(1)}x_{2}^{(2)}%
\end{array}
\]
and
\[%
\begin{array}
[c]{ccccl}%
T_{m,p} & : & \ell_{p}^{2^{m-1}}\times\cdots\times\ell_{p}^{2^{m-1}} &
\rightarrow & \mathbb{R}\\
&  & (x^{(1)},...,x^{(m)}) & \mapsto & \left(  x_{1}^{(m)}+x_{2}^{(m)}\right)
T_{m-1,p}\left(  x^{(1)},...,x^{(m)}\right) \\
&  &  &  & +\left(  x_{1}^{(m)}-x_{2}^{(m)}\right)  T_{m-1,p}\left(
B^{2^{m-1}}(x^{(1)}),B^{2^{m-2}}(x^{(2)}),...,B^{2}(x^{(m-1)})\right)  ,
\end{array}
\]
where $x^{(k)}=\left(  x_{j}^{(k)}\right)  _{j=1}^{{2^{m-1}}}\in\ell
_{p}^{2^{m-1}}$, $1\leq k\leq m$, and $B$ is the backward shift operator in
$\ell_{p}^{2^{m-1}}$. Observe that
\begin{align*}
\left\vert T_{m,p}(x^{(1)},...,x^{(m)})\right\vert  &  \leq\left\vert
x_{1}^{(m)}+x_{2}^{(m)}\right\vert \left\vert T_{m-1,p}\left(  x^{(1)}%
,...,x^{(m)}\right)  \right\vert \\
&  +\left\vert x_{1}^{(m)}-x_{2}^{(m)}\right\vert \left\vert T_{m-1,p}\left(
B^{2^{m-1}}(x^{(1)}),B^{2^{m-2}}(x^{(2)}),...,B^{2}(x^{(m-1)})\right)
\right\vert \\
&  \leq\left\Vert T_{m-1,p}\right\Vert \left(  \left\vert x_{1}^{(m)}%
+x_{2}^{(m)}\right\vert +\left\vert x_{1}^{(m)}-x_{2}^{(m)}\right\vert \right)
\\
&  = \left\Vert T_{m-1,p}\right\Vert 2\max\left\{  \left\vert x_{1}%
^{(m)}\right\vert ,\left\vert x_{2}^{(m)}\right\vert \right\} \\
&  \leq2\left\Vert T_{m-1,p}\right\Vert \left\Vert x^{(m)}\right\Vert _{p}.
\end{align*}
Therefore,
\begin{equation}
\Vert T_{m,p}\Vert\leq2^{m-2}\Vert T_{2,p}\Vert. \label{002}%
\end{equation}

\bigskip

\noindent\textbf{Step 2.} Estimating $\Vert T_{2,4}\Vert$.

\bigskip

Note that
\[
\left\Vert T_{2,4}\right\Vert =\sup\left\{  \left\Vert T_{2,4}^{(x^{(1)}%
)}\right\Vert :\left\Vert x^{(1)}\right\Vert _{4}=1\right\}  ,
\]
where $T_{2,4}^{(x^{(1)})}:\ell_{4}^{2}\rightarrow\mathbb{R}$ is given by
$x^{(2)}\mapsto T_{2,4}\left(  x^{(1)},x^{(2)}\right)  $. Thus we have the
operator
\[
T_{2,4}^{(x^{(1)})}\left(  x^{(2)}\right)  =\left(  x_{1}^{(1)}+x_{2}%
^{(1)}\right)  x_{1}^{(2)}+\left(  x_{1}^{(1)}-x_{2}^{(1)}\right)  x_{2}%
^{(2)}.
\]
Since $\left(  \ell_{4}\right)  ^{\ast}=\ell_{\frac{4}{3}}$, we obtain
$\left\Vert T_{2,4}^{(x^{(1)})}\right\Vert =\left\Vert \left(  x_{1}%
^{(1)}+x_{2}^{(1)},x_{1}^{(1)}-x_{2}^{(1)},0,0,...\right)  \right\Vert
_{\frac{4}{3}}$. Therefore
\[
\left\Vert T_{2,4}\right\Vert =\sup\left\{  \left(  \left\vert x_{1}%
^{(1)}+x_{2}^{(1)}\right\vert ^{\frac{4}{3}}+\left\vert x_{1}^{(1)}%
-x_{2}^{(1)}\right\vert ^{\frac{4}{3}}\right)  ^{\frac{3}{4}}:\left\vert
x_{1}^{(1)}\right\vert ^{4}+\left\vert x_{2}^{(1)}\right\vert ^{4}=1\right\}
.
\]
We can verify that it is enough to maximize the above expression when
$x_{1}^{(1)},x_{2}^{(1)}\geq0$. Then
\begin{align*}
\left\Vert T_{2,4}\right\Vert  &  =\sup\left\{  \left(  \left(  x+\left(
1-x^{4}\right)  ^{\frac{1}{4}}\right)  ^{\frac{4}{3}}+\left\vert x-\left(
1-x^{4}\right)  ^{\frac{1}{4}}\right\vert ^{\frac{4}{3}}\right)  ^{\frac{3}%
{4}}:x\in\left[  0,1\right]  \right\} \\
&  =\max\left\{  \sup\left\{  f(x):x\in\left[  0,2^{-\frac{1}{4}}\right]
\right\}  ,\sup\left\{  g(x):x\in\left[  2^{-\frac{1}{4}},1\right]  \right\}
\right\}
\end{align*}
where
\[
f\left(  x\right)  :=\left(  \left(  x+\left(  1-x^{4}\right)  ^{\frac{1}{4}%
}\right)  ^{\frac{4}{3}}+\left(  \left(  1-x^{4}\right)  ^{\frac{1}{4}%
}-x\right)  ^{\frac{4}{3}}\right)  ^{\frac{3}{4}}%
\]
and
\[
g\left(  x\right)  :=\left(  \left(  x+\left(  1-x^{4}\right)  ^{\frac{1}{4}%
}\right)  ^{\frac{4}{3}}+\left(  x-\left(  1-x^{4}\right)  ^{\frac{1}{4}%
}\right)  ^{\frac{4}{3}}\right)  ^{\frac{3}{4}}.
\]
Examining the maps $f$ and $g$ we easily conclude that%
\[
\left\Vert T_{2,4}\right\Vert <1.74.
\]
In fact, the precise value seems to be graphically $\sqrt{3}$.

\begin{figure}[h]
\centering
\includegraphics[scale=0.25]{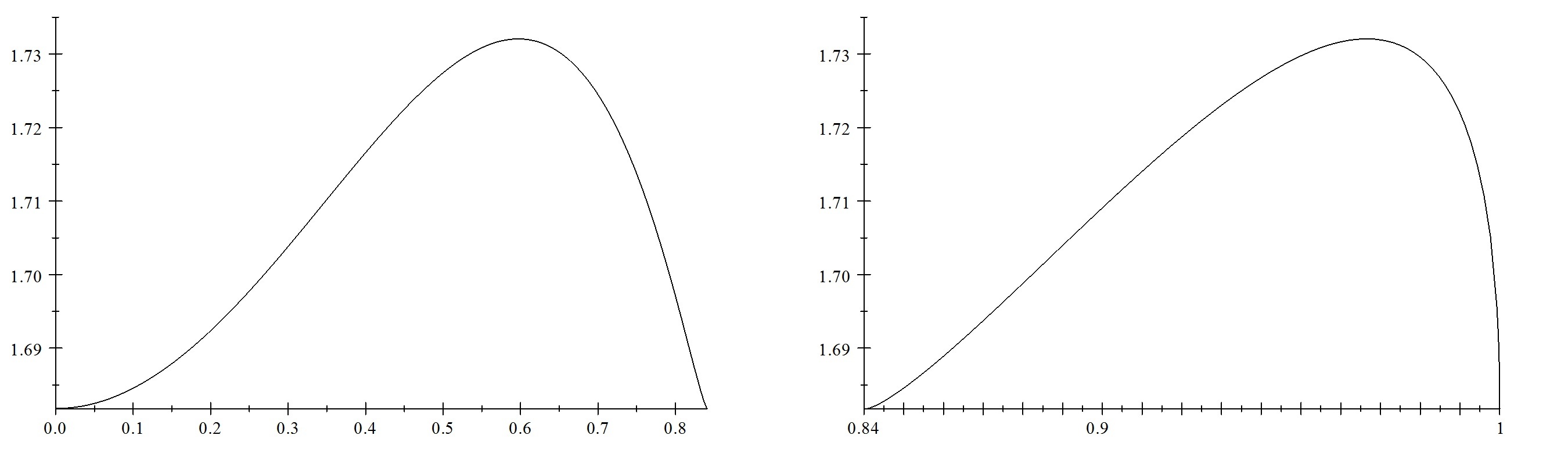}\caption{Graphs of the functions $f$
and $g$, respectively.}%
\label{figure1}%
\end{figure}

\bigskip

\noindent\textbf{Step 3. }Estimating\textbf{ }$\Vert T_{2,p}\Vert$ for
$p\geq4$.

\bigskip

For $1\leq p\leq+\infty$ and $A:\ell_{p}^{N}\times\ell_{p}^{N}\rightarrow
\mathbb{K}$, consider $\widetilde{A}:\ell_{p}^{N}\rightarrow(\ell_{p}%
^{N})^{\ast}$ given by $\widetilde{A}(x)=A(x,\cdot)\in(\ell_{p}^{N})^{\ast}$.
It is easy to see that $\Vert A\Vert=\Vert\widetilde{A}\Vert.$

From \cite{diniz} we know that%
\[
\left\Vert T_{2,\infty}\right\Vert =2.
\]
Let us suppose, for a moment, that we have, in this case, the Riesz-Thorin
Theorem for real scalars with constant $1$, as in the case of complex scalars. By considering
\[
\theta=\frac{p-4}{p},
\]
we would conclude from \cite[Theorem 1.1.1]{bergh} that
\[
\left\Vert \widetilde{T_{2,p}}\right\Vert \leq\left\Vert \widetilde{T_{2,4}%
}\right\Vert ^{1-\theta}\left\Vert \widetilde{T_{2,\infty}}\right\Vert
^{\theta},
\]
i.e.,
\begin{equation}
\Vert T_{2,p}\Vert\leq\Vert T_{2,4}\Vert^{1-\theta}\Vert T_{2,\infty}%
\Vert^{\theta}<(1.74)^{\frac{4}{p}}2^{\frac{p-4}{p}}.\label{mam}%
\end{equation}

\bigskip

\noindent\textbf{Step 4.} Estimating the constants.

\bigskip

From \eqref{002} and \eqref{mam} we would conclude that $\Vert
T_{m,p}\Vert<2^{m-2}(1.74)^{\frac{4}{p}}2^{\frac{p-4}{p}}$. On the other hand,
from Theorem \ref{000} we have
\[
(4^{m-1})^{\frac{mp+p-2m}{2mp}}=\left(  \sum_{j_{1},...,j_{m}=1}^{2^{m-1}%
}\left\vert T_{m,p}(e_{j_{1}},...,e_{j_{m}})\right\vert ^{\frac{2mp}{mp+p-2m}%
}\right)  ^{\frac{mp+p-2m}{2mp}}<C_{m,p}^{\mathbb{R}}2^{m-2}(1.74)^{\frac
{4}{p}}2^{\frac{p-4}{p}}.
\]
and thus
\[
C_{m,p}^{\mathbb{R}}>\frac{(4^{m-1})^{\frac{mp+p-2m}{2mp}}}{2^{m-2}%
(1.74)^{\frac{4}{p}}2^{\frac{p-4}{p}}}=2^{\frac{mp+(6-4\log_{2}(1.74))m-2m^{2}%
-p}{mp}}.
\]

\bigskip

\noindent\textbf{Step 5.} Verifying that $2^{\frac{mp+(6-4\log_{2}%
(1.74))m-2m^{2}-p}{mp}}>1$.

\bigskip

Indeed,
\[
2^{\frac{mp+(6-4\log_{2}(1.74))m-2m^{2}-p}{mp}}=\frac{2^{\frac{mp+6m-2m^{2}%
-p}{mp}}}{(1.74)^{\frac{4}{p}}}>2^{\frac{mp+2m-2m^{2}-p}{mp}}\geq1.
\]
%This calculation also shows that our result improves %\eqref{001}.

However, if $p=2m$ it remains to consider the case in which the Riesz-Thorin Theorem
holds with a constant bigger than $1;$ for this particular case, we may repeat
the Step 2 for other values of $p.$ We just need to observe that $\left\Vert
T_{2,p}\right\Vert <2$ for $4\leq p<\infty$. In fact, in this case,
\[
C_{m,p}^{\mathbb{R}}\geq\frac{(4^{m-1})^{\frac{mp+p-2m}{2mp}}}{2^{m-2}%
\left\Vert T_{m,p}\right\Vert }>2^{\frac{mp+2m-2m^{2}-p}{mp}}\geq1.
\]

%\bigskip

\bigskip

\noindent\textbf{Acknowledgement:} The authors thank the referees for
important suggestions that improved the final version of this paper.


\begin{thebibliography}{9}                                                                                                %
%\bibitem {alencarmatos} R. Alencar and  M. C. Matos, Some classes of multilinear mappings between Banach spaces, Publ. Dep. Análisis Matematico, Universidad Complutense de Madrid, Section 1, Number \textbf{12}, 1989.


%\bibitem{arr} R. Aron, V. I. Gurariy and J. B. Seoane, Lineability and spaceability of sets of functions on $\mathbb{R}$, Proc. Amer. Math. Soc. 133 (2005), no. 3, 795--803.


%\bibitem {aromlacruz} R. M. Aron, M. Lacruz, R. A. Ryan and A. M. Tonge, The generalized Rademacher functions. Note Mat. \textbf{12} (1992), 15--25.


%\bibitem {jfa}N. Albuquerque, F. Bayart, D. Pellegrino and J.
%Seoane-Sep\'{u}lveda, Sharp generalizations of the multilinear
%Bohnenblust--Hille inequality, J. Funct. Anal. 266 (2014), 3726-3740


%\bibitem {n}N. Albuquerque, F. Bayart, D. Pellegrino and J. Seoane--Sep\'{u}lveda, Optimal Hardy--Littlewood type inequalities for polynomials and multilinear operators, arXiv:1311.3177v3 [math.FA].


%\bibitem {nacib}N. Albuquerque, D. Nunez-Alarc\'{o}n, J. Santos, D.M.
%Serrano-Rodr\'{\i}guez, Absolutely summing multilinear operators via
%interpolation, arXiv:1404.4949 [math.FA] 19 April 2014.


\bibitem {aps}G. Ara\'{u}jo, D. Pellegrino, and D.D.P. Silva, On the upper
bounds for the constants of the Hardy--Littlewood inequality, J. Funct. Anal.
\textbf{267} (2014), 1878--1888.

\bibitem {bohr}F. Bayart, D. Pellegrino, and J.B. Seoane-Sep\'{u}lveda, The
Bohr radius of the $n$--dimensional polydisk is equivalent to $\sqrt{(\log
n)/n}$, Advances in Math. \textbf{264} (2014) 726--746.

\bibitem {bergh}J. Bergh and J. L\"{o}fstr\"{om}, Interpolation spaces: an
introduction, Grundlehren der Mathematischen Wissenschaften \textbf{223},
Berlin-New York: Springer-Verlag, ISBN 3-540-07875-4.

%\bibitem {botpell}G. Botelho, D. Pellegrino, Absolutely summing linear
%operators into spaces with no finite cotype, Bull. Belg. Math. Soc. Simon
%Stevin 16 (2009), 373-378.


%\bibitem {df}A. Defant, D. Popa and U. Schwarting, Coordenatewise multiple
%summing operators on Banach spaces, J. Funct. Anal. \textbf{259} (2010), 220--242.


%\bibitem {Di}J. Diestel, H. Jarchow and A. Tonge, Absolutely summing
%operators, Cambridge University Press, Cambridge, 1995.


%\bibitem {a3} L. Bernal--Gonz\'{a}lez and M. Ordonez-Cabrera, Lineability criteria with applications, J. Funct. Anal., \textbf{266} (2014), 3997--4025.


%\bibitem {bernal}L. Bernal--Gonz\'{a}lez, D. Pellegrino and J. Seoane--Sep\'{u}lveda, Linear subsets of nonlinear sets in topological vector
%spaces, Bull. Amer. Math. Soc. \textbf{51} (2014), 71--130.


%\bibitem {bh}H. F. Bohnenblust and E. Hille, On the absolute convergence of
%Dirichlet series, Ann. of Math. \textbf{32} (1931), 600--622.


%\bibitem {ir}G. Botelho, Cotype and absolutely summing multilinear mappings
%and homogeneous polynomials, Proc. Roy.\ Irish Acad. Sect. A \textbf{97}
%(1997), 145--153.


%\bibitem {studia}G. Botelho, D. Cariello, V. Favaro, D. Pellegrino and J. B. Seoane--Sep\'{u}lveda, Distinguished subspaces of $L_p$ of maximal dimension, Studia Math., \textbf{215} (3) (2013), 261-280.


%\bibitem {REMC2010}G. Botelho, C. Michels and D. Pellegrino, Complex
%interpolation and summability properties of multilinear operators, Rev. Mat.
%Complut. \textbf{23} (2010), 139--161.


%\bibitem {b222} G. Botelho and D. Pellegrino, Absolutely summing polynomials on Banach spaces with unconditional basis. J. Math. Anal. Appl. \textbf{321} (2006), nº. 1, 50--58.


%\bibitem {portmath} G. Botelho and D. Pellegrino, Coincidence situations for absolutely summing non--linear mappings. Port. Math. (N.S.) \textbf{64} (2007), nº. 2, 175--191.


%\bibitem {nach}G. Botelho, D. Pellegrino, When every multilinear mapping is
%multiple summing, Math. Nachr. \textbf{282} (2009), nº. 10, 1414--1422.


%\bibitem {indag} G. Botelho, D. Pellegrino, P. Rueda, Summability and estimates for polynomials and multilinear mappings. Indag. Math. (N.S.) \textbf{19} (2008), nº. 1, 23--31.


%\bibitem {rue}G. Botelho, D. Pellegrino and P. Rueda, Dominated polynomials on
%infinite dimensional spaces. Proc. Amer. Math. Soc. \textbf{138} (2010), n. 1, 209--216.


%\bibitem {rue2}G. Botelho, D. Pellegrino and P. Rueda, Dominated bilinear
%forms and $2$--homogeneous polynomials. Publ. Res. Inst. Math. Sci. \textbf{46} (2010),
%nº. 1, 201--208.


%\bibitem {a2} D. Cariello and J. B. Seoane--Sep\'{u}lveda, Basic sequences and spaceability in $\ell_p $ spaces, J. Funct. Anal., \textbf{266} (2014), 3794--3814.


%\bibitem {Di}J. Diestel, H. Jarchow and A. Tonge, Absolutely summing
%operators, Cambridge University Press, Cambridge, 1995.


%\bibitem {dimant}V. Dimant and P. Sevilla--Peris, Summation of coefficients of polynomials on $\ell_p$ spaces, arXiv:1309.6063v1 [math.FA].


\bibitem {diniz}D. Diniz, G. A. Mu\~{n}oz-Fern\'{a}dez, D. Pellegrino, and J.
B. Seoane-Sep\'{u}lveda, Lower Bounds for the constants in the
Bohnenblust-Hille inequality: the case of real scalars, Proc. Amer. Math. Soc.
\textbf{142} (2014), n. 2, 575--580.

\bibitem {hardy}G. Hardy and J. E. Littlewood, Bilinear forms bounded in space
$[p, q]$, Quart. J. Math. \textbf{5} (1934).

%\bibitem {matos}M. C. Matos, Fully absolutely summing mappings and Hilbert
%Schmidt operators, Collect. Math. \textbf{54} (2003) 111--136.


\bibitem {montanaro}A. Montanaro, Some applications of hypercontractive
inequalities in quantum information theory, J. Math. Phys. \textbf{53} (2012).

\bibitem {nunez}D. N\'{u}\~{n}ez-Alarc\'{o}n, A note on the polynomial
Bohnenblust--Hille inequality, J. Math. Anal. Appl., \textbf{407} (2013) 179--181.

\bibitem {npss}D. N\'{u}\~{n}ez-Alarc\'{o}n, D. Pellegrino, D. M.
Serrano-Rodr\'{\i}guez, and J. B. Seoane--Sep\'{u}lveda, There exist
multilinear Bohnenblust--Hille constants $\left(  C_{n}\right)  _{n=1}%
^{\infty}$ with $\lim_{n\to\infty}\left(  C_{n+1}-C_{n}\right)  =0$, J. Funct.
Anal., \textbf{264} (2013) 429--463.

%\bibitem {Da}D. P\'{e}rez-Garc\'{\i}a, Operadores multilineales absolutamente
%sumantes, Ph.D. thesis, Universidad Complutense de Madrid, 2003.


%\bibitem {pv}D. P\'{e}rez-Garc\'{\i}a and I. Villanueva, Multiple summing
%operators on Banach spaces, J. Math. Anal. Appl. \textbf{285} (2003), 86--96.


\bibitem {pra}T. Praciano--Pereira, On bounded multilinear forms on a class of
$\ell_{p}$ spaces, J. Math. Anal. Appl. \textbf{81} (1981), n. 2, 561--568.

%\bibitem {Marcela1}M. L. V. Souza, Aplica\c{c}\~{o}es multilineares
%completamente absolutamente somantes, PhD Thesis, Universidade Estadual de
%Campinas (UNICAMP), 2003.


%\bibitem {drew}L. Drewnowski, Quasicomplements in $F$--spaces, Studia Math. \textbf{77}
%(1984) 373--391.


%\bibitem {LP}J. Lindenstraus and A. Pe\l czy\'{n}ski, Absolutely summing
%operators in $\mathcal{L}_{p}$--spaces and their applications, Studia Math.
%\textbf{29} (1968), 275--324.


%\bibitem {LT}J. Lindenstraus and L. Tzafriri, Classical Banach Spaces I and
%II, Springer--Verlag, 1977.


%\bibitem {LLL}J. E. Littlewood, On bounded bilinear forms in an infinite number
%of variables, Quart. J. (Oxford Ser.) \textbf{1} (1930), 164--174.


%\bibitem {KT} D. Kitson and R. M. Timoney, Operator ranges and spaceability, J. Math. Anal. Appl. \textbf{378}, 2 (2011) 680--686.


%\bibitem {matos2} M. C. Matos, On multilinear mappings of nuclear type. Rev. Mat. Univ. Complut. Madrid \textbf{6} (1993), nº. 1, 61--81.


%\bibitem {mi}B.S. Mitiagin and A. Pe\l czy\'{n}ski, Nuclear operators and
%approximative dimension. Proc. Inter. Congr. of Math. Moscow 1966


%\bibitem {St}D. Pellegrino, Cotype and absolutely summing homogeneous
%polynomials in $\mathcal{L}_{p}$--spaces, Studia Math. \textbf{157} (2003) 121--131.


%\bibitem {laa} D. Pellegrino, Sharp coincidences for absolutely summing
%multilinear operators, Linear Algebra and its Applications, \textbf{440} (2014) 188--196.


%\bibitem {comparing}D. P\'{e}rez--Garc\'{\i}a, Comparing different classes of
%absolutely summing multilinear operators, Arch. Math. (Basel) \textbf{85}
%(2005), 258--267.


%\bibitem {per}D. P\'{e}rez--Garc\'{\i}a, Operadores multilineales absolutamente sumantes, Thesis, Universidad Complutense de Madrid, 2003.


%\bibitem {arc}D. P\'{e}rez--Garc\'{\i}a and I. Villanueva, Multiple summing
%operators on $C(K)$ spaces, Ark. Mat. \textbf{42} (2004), 153--171.


%\bibitem {piet}A. Pietsch, Absolut $p$--summierende Abbildungen in normieten
%Raumen, Studia Math. \textbf{27} (1967), 333--353.


%\bibitem {pi1}A. Pietsch, Ideals of multilinear functionals, Proceedings of
%the Second International Conference on Operator Algebras, Ideals and Their
%Applications in Theoretical Physics, Teubner--texte Math. \textbf{67} (Teubner, Leipzig,
%1983) 185--199.

\end{thebibliography}
\end{document}